\documentclass{article}
\usepackage{amssymb}
\usepackage{amsmath}

\input{tcilatex}
\begin{document}

\begin{center}
%Folder: 410PropZerosWilsonPol; File: 410PropZerosWilsonPolLettersMPSubmission

\bigskip

{\LARGE Properties of the zeros of the polynomials belonging to the Askey
scheme}

\bigskip

$^{\ast }$\textbf{Oksana Bihun}$^{1}$ and $^{+\lozenge }$\textbf{Francesco
Calogero}$^{2}\bigskip $

$^{\ast }$Department of Mathematics, Concordia College at Moorhead,\\
901 8th Str South, Moorhead, MN 56562, USA, +1-218-299-4396

\smallskip

$^{+}$Physics Department, University of Rome \textquotedblleft La Sapienza",\\ p. Aldo Moro, I-00185 ROMA, Italy, +39-06-4991-4372

$^{\lozenge }$Istituto Nazionale di Fisica Nucleare, Sezione di Roma
\smallskip

$^{1}$obihun@cord.edu
\smallskip

$^{2}$francesco.calogero@roma1.infn.it, francesco.calogero@uniroma1.it

\bigskip

\textit{Abstract}
\end{center}

In this paper we provide properties---which are, to the best of our
knowledge, new---of the zeros of the polynomials belonging to the Askey
scheme. These findings include Diophantine relations satisfied by these
zeros when the parameters characterizing these polynomials are appropriately
restricted.

\smallskip


\smallskip

MSC 33C45,
%Hypergeometric functions, Orthogonal polynomials, Askey scheme
11D41, 
%Higher degree Diophantine equations; Fermat's equation
15A18 
%Linear algebra,  	Eigenvalues, singular values, and eigenvectors
\bigskip

\section{Introduction}

The properties of the zeros of polynomials are a core problem of mathematics
to which, over time, an immense number of investigations have been devoted.
Nevertheless new findings in this area continue to emerge, see, for
instance, \cite{BCD1,BCD2,BCD3,BCD4,BCD5,CI2013,IR2013,BC2013,BC2014,CY1, CY2,BCY,C1,C2,C3}. In this paper we
report findings concerning the zeros of the polynomials belonging to the
Askey scheme (see for instance \cite{KS}). These results are, to the best of
our knowledge, \textit{new}; although the main approach to uncover them is
not new, except for some not quite trivial variations. These findings
include Diophantine relations satisfied by these zeros when the parameters
characterizing these polynomials are appropriately restricted.

The technique used to arrive at the results reported below may be considered
to originate from the possibility---firstly noted by G. Sz\"{e}go---to
relate the zeros of classical polynomials to the equilibria of certain
many-body problems \cite{S1939}, and from the subsequent observation that
the investigation of the zeros of \textit{time-dependent} polynomials
satisfying \textit{linear} PDEs provides a method to identify \textit{%
solvable} many-body problems (see \cite{C1978}, and for related developments
for instance the two books \cite{C2001,C2008}; and see for instance \cite%
{ABCOP1979} for many findings which originated from this development and are
analogous to those reported below but concern mainly the zeros of \textit{%
classical} polynomials).

The main findings of this paper are reported in the following Section 2.
They detail properties of the zeros of the Wilson and Racah polynomials,
which are the two \textquotedblleft highest" classes of polynomials
belonging to the Askey scheme \cite{KS}---so that these polynomials feature $%
4$ arbitrary parameters in addition to their degree $N$ (an arbitrary
positive integer). Analogous properties can of course be obtained, from the
results reported below, for the zeros of the (variously \textquotedblleft
named" \cite{KS}) polynomials belonging to \textquotedblleft lower" classes
of the Askey scheme, via the reductions---corresponding to special
assignments of the parameters---that characterize the Askey scheme \cite{KS}%
; but we leave this as a task for the interested reader. The results thus
obtained for the \textquotedblleft lowest" classes of polynomials belonging
to the Askey scheme shall of course reproduce already known findings valid
for the \textit{classical} polynomials \cite{ABCOP1979}.

Our findings are proven in Section 3, and a very terse Section 4 entitled
Outlook outlines future developments. The definitions and some standard
properties of the Wilson and Racah polynomials are reported in the Appendix,
for the convenience of the reader and also to specify our notation; the
reader is advised to glance through this Appendix before reading the next
section, and then to return to it whenever appropriate.

The results that follow are only a consequence of the \textit{explicit} 
\textit{definitions} of the Wilson and Racah polynomials and of the \textit{%
difference equations} they satisfy (see the Appendix); the orthogonality
properties that these polynomials satisfy---in the case of the Racah
polynomials, for integrations on a discrete measure---play no role, so that
the results reported below do \textit{not} require the restrictions on the
parameters of these polynomials---or, in the case of the Racah polynomials,
on their arguments---that are instead mandatory for the validity of their
orthogonality properties.

Finally, let us emphasize that the results reported in this paper are 
\textit{not} special cases of those reported in \cite{BC2014}; indeed the
Wilson and Racah polynomials, while defined in terms of generalized
hypergeometric functions, do not belong to the class of generalized
hypergeometric polynomials treated in \cite{BC2014}, due to the quite
different definitions of the quantity playing the role of argument of these
polynomials.

\section{Main results}

To formulate our main results we refer to the definitions and standard
properties of the Wilson and Racah polynomials as reported in the Appendix,
to which the reader should also refer for the notation employed hereafter.
The proofs of these results are provided in the following Section 3.

Hereafter the symbol $\mathbf{i}$ is the \textit{imaginary unit}, $\mathbf{i}%
^{2}=-1.$

Our main results concerning the zeros of Wilson and Racah polynomials read
as follows:

%\begin{proposition}
%\label{prop2pt1} 

\vspace{2mm}

\noindent \textit{Proposition 2.1}. Let $\bar{z}_{s}=\bar{x}_{s}^{2}$ (with $%
s=1,2,...,N)$ be the $N$ zeros of the Wilson polynomial $W_{N}\left(
z;a,b,c,d\right) \equiv W_{N}\left( x^{2};a,b,c,d\right) $ of degree $N$ in $%
z=x^{2}$ (with $N$ an arbitrary positive integer), so that $W_{N}\left( \bar{%
z}_{s};a,b,c,d\right) \equiv W_{N}\left( \bar{x}_{s}^{2};a,b,c,d\right) =0$
for $s=1,2,...,N$. Define the $N\times N$ matrix $\underline{M}$,
componentwise, as follows: 
\begin{subequations}
\label{Mnm}
\begin{eqnarray}
M_{nn} &=&\left( 2\bar{x}_{n}\right) ^{-2}\Bigg\{\left[ \frac{2A\left( \bar{x%
}_{n}\right) }{\mathbf{i}x_{n}}+\mathbf{i}A^{\prime }\left( \bar{x}%
_{n}\right) \right] ~\dprod\limits_{\ell =1,~\ell \neq n}^{N}\left( 1-\frac{%
1+2\mathbf{i}\bar{x}_{n}}{\bar{x}_{n}^{2}-\bar{x}_{\ell }^{2}}\right)  \notag
\\
&&+2A\left( \bar{x}_{n}\right) ~\sum_{m=1,~m\neq n}\frac{\mathbf{i}\bar{x}%
_{n}-\left( \bar{x}_{n}^{2}+\bar{x}_{m}^{2}\right) }{\left( \bar{x}_{n}^{2}-%
\bar{x}_{m}^{2}\right) ^{2}}\dprod\limits_{\ell =1,~\ell \neq n,m}^{N}\left(
1-\frac{1+2\mathbf{i}\bar{x}_{n}}{\bar{x}_{n}^{2}-\bar{x}_{\ell }^{2}}\right)
\notag \\
&&~+\left[ \left( \bar{x}_{s}\rightarrow (-\bar{x}_{s})\right) \right] ~%
\Bigg\}~,~~~n=1,2,...,N~,
\end{eqnarray}%
\begin{eqnarray}
&&M_{nm}=-\left( 2\bar{x}_{n}\right) ^{-2}\Bigg\{2A\left( \bar{x}_{n}\right)
~\frac{\mathbf{i}\bar{x}_{m}\left( 1+2\mathbf{i}\bar{x}_{n}\right) }{\left( 
\bar{x}_{n}^{2}-\bar{x}_{m}^{2}\right) ^{2}}\dprod\limits_{\ell =1,~\ell
\neq n,m}^{N}\left( 1-\frac{1+2\mathbf{i}\bar{x}_{n}}{\bar{x}_{n}^{2}-\bar{x}%
_{\ell }^{2}}\right)  \notag \\
&&~+\left[ \left( \bar{x}_{s}\rightarrow (-\bar{x}_{s})\right) \right] ~%
\Bigg\}\,,~~~m,n=1,2,...,N,~~~m\neq n~,
\end{eqnarray}%
where 
\end{subequations}
\begin{subequations}
\label{AAprime}
\begin{equation}
A\left( \bar{x}_{n}\right) \equiv A\left( \bar{x}_{n};a,b,c,d\right) =\alpha
_{4}+\mathbf{i}\alpha _{3}\bar{x}_{n}-\alpha _{2}\bar{x}_{n}^{2}-\mathbf{i}%
\alpha _{1}\bar{x}_{n}^{3}+\bar{x}_{n}^{4}~,  \label{A}
\end{equation}%
\begin{equation}
A^{\prime }\left( \bar{x}_{n}\right) \equiv A\left( \bar{x}%
_{n};a,b,c,d\right) =\mathbf{i}\alpha _{3}-2\alpha _{2}\bar{x}_{n}-3\mathbf{i%
}\alpha _{1}\bar{x}_{n}^{2}+4\bar{x}_{n}^{3}~,
\end{equation}%
\end{subequations}
\begin{subequations}
\label{alpha}
\begin{equation}
\alpha _{1}\equiv \alpha _{1}\left( a,b,c,d\right) =a+b+c+d~,  \label{alpha1}
\end{equation}%
\begin{equation}
\alpha _{2}\equiv \alpha _{2}\left( a,b,c,d\right) =ab+ac+ad+bc+bd+cd~,
\end{equation}%
\begin{equation}
\alpha _{3}\equiv \alpha _{3}\left( a,b,c,d\right) =bcd+acd+abd+abc~,
\end{equation}%
\begin{equation}
\alpha _{4}\equiv \alpha _{4}\left( a,b,c,d\right) =abcd~.
\end{equation}%
The symbol $+\left[ \left( \bar{x}_{s}\rightarrow (-\bar{x}_{s})\right) %
\right] $ denotes the \textit{addition} of everything that comes before it
(within the curly brackets), with the replacement of $\bar{x}_{s}$ with $(-%
\bar{x}_{s})$ for all $s=1,2,\ldots ,N$. Note that these matrix elements $%
M_{nm}$ are \textit{even} functions of the $\bar{x}_{s}$ (i. e., functions
of $\bar{z}_{s}=\bar{x}_{s}^{2}$; rather than $\bar{x}_{s}$), because the
addition implied by the symbol $+\left[ \left( \bar{x}_{s}\rightarrow (-\bar{%
x}_{s})\right) \right] $ causes \textit{all} terms \textit{odd} in $\bar{x}%
_{s}$ to cancel out (for all values of $s=1,2,...,N$).

Then this $N\times N$ matrix $\underline{M}$ has the $N$ eigenvalues 
\end{subequations}
\begin{equation}
m\left( 2N-m+a+b+c+d-1\right) =m\left( 2N-m+\alpha _{1}-1\right)
,\;\;\;m=1,2,...,N~.~\square  \label{eq:eigM}
\end{equation}
%\end{proposition}

%\begin{proposition}
%. \label{prop2pt2} 

\vspace{2mm}

\noindent \textit{Proposition 2.2.} Let $R_{N}\left( z;\alpha ,\beta ,\gamma
,\delta \right) \equiv R_{N}\left( y^{2}-\theta ^{2};\alpha ,\beta ,\gamma
,\delta \right) $ be the Racah polynomial of degree $N\leq \nu $ in the
variable $z=y^{2}-\theta ^{2}$, where $\theta =\frac{\gamma +\delta +1}{2}$,
with the parameters $\alpha ,\beta ,\gamma ,\delta $ and $\nu $ satisfying
conditions \eqref{condRacah}. Let $\bar{z}_{s}=\bar{y}_{s}^{2}-\theta ^{2}$
be the $N$ zeros of this polynomial, so that $R_{N}\left( \bar{z}_{s};\alpha
,\beta ,\gamma ,\delta \right) =R_{N}\left( \bar{y}_{s}^{2}-\theta
^{2};\alpha ,\beta ,\gamma ,\delta \right) =0$ for all $s=1,2,...,N$. Define
the $N\times N$ matrix $\underline{\tilde{M}}$, componentwise, as follows: 
\begin{subequations}
\label{tildeM}
\begin{eqnarray}
&&\tilde{M}_{nn}=\frac{1}{2}\Bigg\{\left[ \left( \frac{\tilde{D}(\bar{y}_{n})%
}{\bar{y}_{n}^{2}}-\frac{\tilde{D}^{\prime }(\bar{y}_{n})}{\bar{y}_{n}}%
\right) (1+2\bar{y}_{n})-2\frac{\tilde{D}(\bar{y}_{n})}{\bar{y}_{n}}\right]
\prod_{\ell =1,\ell \neq n}^{N}\left( 1+\frac{1+2\bar{y}_{n}}{\bar{y}%
_{n}^{2}-\bar{y}_{\ell }^{2}}\right)   \notag \\
&&+2\frac{\tilde{D}(\bar{y}_{n})}{\bar{y}_{n}}(1+2\bar{y}_{n})\sum_{m=1,m%
\neq n}^{N}\left[ \frac{\bar{y}_{n}^{2}+\bar{y}_{m}^{2}+\bar{y}_{n}}{(\bar{y}%
_{n}^{2}-\bar{y}_{m}^{2})^{2}}\prod_{\ell =1,\ell \neq n,m}^{N}\left( 1+%
\frac{1+2\bar{y}_{n}}{\bar{y}_{n}^{2}-\bar{y}_{\ell }^{2}}\right) \right]  
\notag \\
&&+[(\bar{y}_{s}\rightarrow (-\bar{y}_{s}))]\Bigg\},\;\;n=1,2,\ldots ,N
\label{tildeMnn}
\end{eqnarray}%
and 
\begin{eqnarray}
\tilde{M}_{nm} &=&-\frac{1}{(\bar{y}_{n}^{2}-\bar{y}_{m}^{2})^{2}}\Bigg\{%
\frac{\bar{y}_{m}\tilde{D}(\bar{y}_{n})}{\bar{y}_{n}}(1+2\bar{y}%
_{n})^{2}\prod_{\ell =1,\ell \neq n,m}^{N}\left( 1+\frac{1+2\bar{y}_{n}}{%
\bar{y}_{n}^{2}-\bar{y}_{\ell }^{2}}\right)   \notag \\
&&+[(\bar{y}_{s}\rightarrow (-\bar{y}_{s}))]\Bigg\},\;\;n,m=1,2,\ldots
,N,~~~n\neq m,  \label{tildeMnm}
\end{eqnarray}%
where 
\end{subequations}
\begin{equation}
\tilde{D}(\bar{y}_{n})=\frac{(2\bar{y}_{n}+\gamma +\delta +1)(2\bar{y}%
_{n}+\gamma -\delta +1)(2\bar{y}_{n}+2\alpha -\gamma -\delta +1)(2\bar{y}%
_{n}+2\beta -\gamma +\delta +1)}{32\bar{y}_{n}(2\bar{y}_{n}+1)}
\label{eq:R10}
\end{equation}%
and again, the symbol $+[(\bar{y}_{s}\rightarrow (-\bar{y}_{s}))]$ denotes
the addition of everything that comes before it (within the curly brackets),
with $\bar{y}_{s}$ replaced by $(-\bar{y}_{s})$ for all $s=1,2,\ldots ,N$.
Note that these matrix elements $\tilde{M}_{nm}$ are functions \hspace{0in}%
of $\bar{z}_{s}=\bar{y}_{s}^{2}-\theta ^{2}$, $s=1,2,\ldots ,N$, because the
addition implied by the symbol $+\left[ \left( \bar{y}_{s}\rightarrow (-\bar{%
y}_{s})\right) \right] $ causes \textit{all} terms \textit{odd} in $\bar{y}%
_{s}$ to cancel out (for all values of $s=1,2,...,N$).

Then this $N\times N$ matrix $\underline{\tilde{M}}$ has the eigenvalues 
\begin{equation}
m(m-2N-\alpha -\beta -1),\;\;\;m=1,2,\ldots ,N~.~\square
\label{eq:eigMtilde}
\end{equation}
%\end{proposition}

The following immediate corollaries of \textit{Propositions}~2.1. and~2.2
are worth a mention.

%\begin{corollary}
%. \label{cor2pt1pt1} 
\vspace{2mm}

\noindent \textit{Corollary 2.3}.

%\begin{itemize}
%\item[(i)] 
\noindent (i) (Wilson) If $\alpha _{1}=a+b+c+d$ (see (\ref{alpha1})) is an 
\textit{integer} number (or a \textit{rational} number), the $N$ eigenvalues
of the $N\times N$ matrix $\underline{M}$ (see (\ref{Mnm})) are all \textit{%
integer }(or \textit{rational}) numbers.

%\item[(ii)] 
\noindent (ii) (Racah) If $\alpha +\beta $ is an \textit{integer} number (or
a \textit{rational} number), the $N$ eigenvalues of the $N\times N$ matrix $%
\underline{\tilde{M}}$ (see (\ref{tildeM})) are all \textit{integer }(or 
\textit{rational}) numbers. $\square $. %\end{itemize}
%\end{corollary}

These are remarkable Diophantine properties. \pagebreak

%\begin{corollary}
%. \label{cor2pt1pt2} 

\vspace{2mm}

\noindent \textit{Corollary 2.4}.

%\begin{itemize}
%\item[(i)] 
\noindent (i) (Wilson) The $N\times N$ matrix $\underline{M}$ (see (\ref{Mnm}%
))---which depends of course on the $4$ \textit{a priori} arbitrary
parameters $a,b,c,d$, explicitly via $A\left( \bar{x}_{n}\right) $ and $%
A^{\prime }\left( \bar{x}_{n}\right) ,$ see (\ref{AAprime}) and (\ref{alpha}%
), and implicitly via the dependence on these $4$ parameters of the $N$
zeros $\bar{z}_{n}\equiv \bar{z}_{n}\left( N;a,b,c,d\right) $ of the Wilson
polynomials $W_{N}\left( z;a,b,c,d\right) \equiv W_{N}\left(
x^{2};a,b,c,d\right) $---is \textit{isospectral} under any variation of
these $4$ parameters which does not change the value of $\alpha _{1}\equiv
a+b+c+d$---or, equivalently, under any variation (without any restriction)
of the values of the $3$ parameters $\alpha _{2},\alpha _{3},\alpha _{4}$
(see (\ref{alpha})).

%\item[(ii)] 
\noindent (ii) (Racah) The $N\times N$ matrix $\underline{\tilde{M}}$
defined by~\eqref{tildeM}, which depends on the $4$ parameters $\alpha
,\beta ,\gamma ,\delta $ satisfying condition~\eqref{condRacah}, explicitly
via $\tilde{D}(\bar{y}_{n})$ and implicitly via the dependence of the zeros $%
\bar{z}_{n}=\bar{y}_{n}^{2}-\theta ^{2}$ of the Racah polynomial $%
R_{N}(z;\alpha ,\beta ,\gamma ,\delta )=R_{N}(y^{2}-\theta ^{2};\alpha
,\beta ,\gamma ,\delta )$ on these $4$ parameters, is \textit{isospectral}
under any variation of these $4$ parameters that does not change the value
of $\alpha +\beta $ and preserves condition \eqref{condRacah}. $\square $ 
%\end{itemize}
%\end{corollary}

%\begin{corollary}
%. \label{cor2pt1pt3} 
\vspace{2mm}

\noindent \textit{Corollary 2.5}. Several identities satisfied by the zeros
of the Wilson and the Racah polynomials are implied by the following
standard consequences of \textit{Propositions}~2.1 and~2.2:

%\begin{itemize}
%\item[(i)] 
\noindent (i) (Wilson) 
\begin{subequations}
\begin{equation}
\func{trace}\left[ \left( \underline{M}\right) ^{k}\right] =\sum_{m=1}^{N}%
\left[ m\left( 2N-m+\alpha _{1}-1\right) \right] ^{k}~,~~~k=1,2,3,...~,
\end{equation}
\begin{equation}
\det \left[ \underline{M}\right] =\dprod\limits_{m=1}^{N}\left[ m\left(
2N-m+\alpha _{1}-1\right) \right] =\frac{N!~\Gamma \left( 2N+\alpha
_{1}-1\right) }{\Gamma \left( N+\alpha _{1}-1\right) }~.
\end{equation}

%\item[(ii)] 
\noindent (ii) (Racah) 
\end{subequations}
\begin{subequations}
\begin{equation}
\func{trace}\left[ \left( \underline{\tilde{M}}\right) ^{k}\right]
=\sum_{m=1}^{N}\left[ m(m-2N-\alpha -\beta -1)\right] ^{k}~,~~~k=1,2,3,...~,
\end{equation}%
\begin{equation}
\det \left[ \underline{\tilde{M}}\right] =\dprod\limits_{m=1}^{N}\left[
m(m-2N-\alpha -\beta -1)\right] =\frac{N!\Gamma (-N-\alpha -\beta )}{\Gamma
(-2N-\alpha -\beta )}~.~\square
\end{equation}
%\end{itemize}

%\end{corollary}

\section{Proofs of the main results}

We begin with the proof of \textit{Proposition~2.1}.

%\begin{proof}
Let $\psi _{2N}\left( x,t\right) \equiv \psi _{2N}\left( x,t;a,b,c,d\right) $
be an \textit{even} \textit{monic polynomial} of degree $2N$ in $x$
satisfying the Differential-Difference Equation (DDE)

\end{subequations}
\begin{subequations}
\label{DDE}
\begin{equation}
\frac{\partial ~\psi _{2N}\left( x,t\right) }{\partial ~t}=\mathbf{i~}\left[
B\left( x\right) \left( 1-\delta _{-}\right) +B\left( -x\right) \left(
1-\delta _{+}\right) +N(N+\alpha _{1}-1)\right] \psi _{2N}\left( x,t\right)
~,  \label{DDEa}
\end{equation}%
where $B\left( x\right) \equiv B\left( x;a,b,c,d\right) $ is defined by~%
\eqref{B+-} and the shift operators $\delta _{\pm }$ act on functions of the
variable $x$ as follows: 
\begin{equation}
\delta _{\pm }~f\left( x\right) =f\left( x\pm \mathbf{i}\right) ~,
\label{delta+-}
\end{equation}%
and $\alpha _{1}=a+b+c+d$, see (\ref{alpha1}).

To prove that such polynomials exist we introduce their representation as a
sum (with $t$-dependent coefficients) over the first $N+1$ (monic) Wilson
polynomials: 
\end{subequations}
\begin{equation}
\psi _{2N}\left( x,t\right) =p_{N}\left( x^{2}\right) +\sum_{m=1}^{N}\left[
c_{m}\left( t\right) ~p_{N-m}\left( x^{2}\right) \right] ~.  \label{psicm}
\end{equation}%
Here of course $p_{\ell }\left( z\right) \equiv p_{\ell }\left(
z;a,b,c,d\right) $ is the \textit{monic} Wilson polynomial of degree $\ell $%
, see (\ref{MonicW}).

It is then plain that this representation is valid---i. e., consistent with
the \textit{monic polynomial} character of the solution $\psi _{2N}\left(
x,t\right) $ of the DDE\ (\ref{DDE})---because its insertion in the DDE\ (%
\ref{DDE}) implies, via DDE (\ref{DiffEq}), the following simple $t$%
-evolution of the $N$ coefficients $c_{m}$: 
\begin{subequations}
\label{EvEqcm}
\begin{equation}
\dot{c}_{m}\left( t\right) =\mathbf{i}~m~\left( 2N-m+\alpha _{1}-1\right)
~c_{m}\left( t\right) ~,
\end{equation}%
entailing%
\begin{equation}
c_{m}\left( t\right) =c_{m}\left( 0\right) ~\exp \left[ \mathbf{i}~m~\left(
2N-m+\alpha _{1}-1\right) ~t\right] ~.  \label{cmt}
\end{equation}%
Of course here and hereafter a superimposed dot denotes a $t$%
-differentiation.

Hereafter we shall occasionally take advantage of the following obvious
notational identities: 
\end{subequations}
\begin{subequations}
\begin{equation}
\psi _{2N}\left( x,t\right) \equiv \Psi _{N}\left( z,t\right) ~\ ~\text{with}%
~~~z=x^{2}~,
\end{equation}%
\begin{equation}
p_{2N}\left( x\right) \equiv P_{N}\left( z\right) ~\ ~\text{with}~~~z=x^{2}~.
\end{equation}

Next, let us introduce the $N$ (of course $t$-dependent) zeros $z_{n}\left(
t\right) $ of the polynomial $\Psi _{N}\left( z,t\right) $ by writing 
\end{subequations}
\begin{subequations}
\begin{equation}
\Psi _{N}\left( z,t\right) =\prod\limits_{n=1}^{N}\left[ z-z_{n}\left(
t\right) \right] ~,
\end{equation}%
and likewise%
\begin{equation}
\psi _{2N}\left( x,t\right) \equiv \Psi _{N}\left( x^{2},t\right)
=\prod\limits_{n=1}^{N}\left[ x^{2}-x_{n}^{2}\left( t\right) \right] ~,
\label{psizeros}
\end{equation}%
implying of course $z_{n}\left( t\right) =x_{n}^{2}\left( t\right) $.

Let us now focus on the $t$-evolution of these $N$ zeros $x_{n}\equiv
x_{n}\left( t\right) $, as implied by the DDE (\ref{DDE}).

The first observation is that the \textit{equilibrium} configuration of
these $N$ zeros is provided by the $N$ zeros $\bar{z}_{n}=\bar{x}_{n}^{2}$
of the Wilson polynomial $W_{N}\left( z\right) =w_{2N}\left( x\right) $ (of
degree $N$ in $z=x^{2}$). Indeed, if the solution $\psi _{2N}(x,t)$ of DDE~(%
\ref{DDEa}) satisfies $\psi _{2N}(x,0)=w_{2N}(x)$, then the right-hand side
of (\ref{DDEa}) vanishes at $t=0$, which implies $\psi _{2N}(x,t)=\psi
_{2N}(x,0)=w_{2N}(x)$ is time-independent. Hence $x_{n}\left( 0\right) =\bar{%
x}_{n}$ implies $x_{n}\left( t\right) =x_{n}\left( 0\right) =\bar{x}_{n}$.
Likewise, this of course implies that, for the $N$ coefficients $c_{m}$, see
(\ref{psicm}), the \textit{equilibrium} configuration is that \textit{they
all vanish}, i. e. $c_{m}\left( t\right) =c_{m}\left( 0\right) =\bar{c}%
_{m}=0 $; this being indeed the unique \textit{equilibrium} configuration
associated to the simple system of $t$-evolution equations (\ref{EvEqcm}).

Next, let us look at the actual $t$- evolution of the zeros $x_{n}\left(
t\right) $ implied by DDE (\ref{DDE}). To this end the following formulas
are useful: 
\end{subequations}
\begin{subequations}
\label{eq10}
\begin{equation}
\frac{\partial ~\psi _{2N}\left( x,t\right) }{\partial ~t}=-2~\psi
_{2N}\left( x,t\right) ~\sum_{m=1}^{N}\dot{x}_{m}\left( t\right)
~x_{m}\left( t\right) ~\left[ x^{2}-x_{m}^{2}\left( t\right) \right] ^{-1}~,
\end{equation}%
which obtains by logarithmic $t$-differentiation of (\ref{psizeros}), and
clearly implies (again, via (\ref{psizeros})), for $x=x_{n}\left( t\right) $%
, the identity%
\begin{eqnarray}
&&\left. \frac{\partial ~\psi _{2N}\left( x,t\right) }{\partial ~t}%
\right\vert _{x=x_{n}\left( t\right) }\equiv \frac{\partial ~\psi
_{2N}\left( x_{n}\left( t\right) ,t\right) }{\partial ~t}  \notag \\
&=&-2~\dot{x}_{n}\left( t\right) ~x_{n}\left( t\right)
~\prod\limits_{m=1,~m\neq n}^{N}\left[ x_{n}^{2}\left( t\right)
-x_{m}^{2}\left( t\right) \right] ~.  \label{derivtpsi}
\end{eqnarray}

As for the right-hand side of DDE (\ref{DDE}), the following chain of
relations are clearly implied by (\ref{delta+-}), (\ref{alpha1}), (\ref%
{psizeros}) and (\ref{B+-}): 
\end{subequations}
\begin{eqnarray}
&&\left. \left[ B\left( x\right) ~\left( 1-\delta _{-}\right) +B\left(
-x\right) ~\left( 1-\delta _{+}\right) +N~(N+\alpha _{1}-1)\right] ~\psi
_{2N}\left( x,t\right) \right\vert _{x=x_{n}\left( t\right) }  \notag \\
&=&-B\left( x_{n}\left( t\right) \right) ~\psi _{2N}\left( x_{n}\left(
t\right) -\mathbf{i},t\right) -B\left( -x_{n}\left( t\right) \right) ~\psi
_{2N}\left( x_{n}\left( t\right) +\mathbf{i},t\right)   \notag \\
&=&\left[ 1+2\mathbf{i}x_{n}\left( t\right) \right] B\left( x_{n}\left(
t\right) \right) \dprod\limits_{m=1,~m\neq n}^{N}\left[ x_{n}^{2}\left(
t\right) -x_{m}^{2}\left( t\right) -1-2\mathbf{i}x_{n}\left( t\right) \right]
\notag \\
&&+\left[ 1-2\mathbf{i}x_{n}\left( t\right) \right] B\left( -x_{n}\left(
t\right) \right) \dprod\limits_{m=1,~m\neq n}^{N}\left[ x_{n}^{2}\left(
t\right) -x_{m}^{2}\left( t\right) -1+2\mathbf{i}x_{n}\left( t\right) \right]
\notag \\
&=&\left\{ \frac{\left( a+\mathbf{i}x_{n}\right) \left( b+\mathbf{i}%
x_{n}\right) \left( c+\mathbf{i}x_{n}\right) \left( d+\mathbf{i}x_{n}\right) 
}{2\mathbf{i}x_{n}}\dprod\limits_{m=1,~m\neq n}^{N}\left[
x_{n}^{2}-x_{m}^{2}-1-2\mathbf{i}x_{n}\right] \right.   \notag \\
&&\left. -\frac{\left( a-\mathbf{i}x_{n}\right) \left( b-\mathbf{i}%
x_{n}\right) \left( c-\mathbf{i}x_{n}\right) \left( d-\mathbf{i}x_{n}\right) 
}{2\mathbf{i}x_{n}}\dprod\limits_{m=1,~m\neq n}^{N}\left[
x_{n}^{2}-x_{m}^{2}-1+2\mathbf{i}x_{n}\right] \right\} ~.  \notag
\label{eq11} \\
&&
\end{eqnarray}%
Note that, in the last right-hand side, we omitted for notational simplicity
to indicate explicitly the $t$-dependence of the quantities $x_{n}\equiv
x_{n}\left( t\right) $ and $x_{m}\equiv x_{m}\left( t\right) $; and we shall
often also do so below.

Hence it is now clear---by inserting this last formula and (\ref{derivtpsi})
in DDE~(\ref{DDE})---that the $t$-dependent zeros $z_{n}\left( t\right)
=x_{n}^{2}\left( t\right) $ of the polynomial $\psi _{2N}\left( x,t\right)
\equiv \Psi _{N}\left( x^{2},t\right) \equiv \Psi _{N}\left( z,t\right) $
satisfy the following system of \textquotedblleft equations of motion": 
\begin{eqnarray}
&&2x_{n}\dot{x}_{n}=\dot{z}_{n}=-\mathbf{i~}\Bigg\{\frac{A\left(
x_{n}\right) }{2\mathbf{i}x_{n}}~\left[ \dprod\limits_{m=1,~m\neq
n}^{N}\left( \frac{x_{n}^{2}-x_{m}^{2}-1-2\mathbf{i}x_{n}}{%
x_{n}^{2}-x_{m}^{2}}\right) \right]   \notag \\
&&+\left[ \left( x_{s}\rightarrow (-x_{s})\right) \right] \Bigg\}  \notag \\
&=&F_{n}\left( x_{n}^{2};x_{\ell }^{2},\ell \neq n\right) =F_{n}\left(
z_{n};z_{\ell },\ell \neq n\right) =F_{n}\left( z_{n};z_{\ell },\ell \neq
n;a,b,c,d\right) ~,  \notag \\
&&  \label{EvEqxn}
\end{eqnarray}%
where the $4$ parameters $\alpha _{1},\alpha _{2},\alpha _{3},\alpha _{4}$
are defined in terms of the $4$ parameters $a,b,c,d$ as above, see (\ref%
{alpha}), the symbol $\left[ \left( x_{s}\rightarrow (-x_{s})\right) \right] 
$ has the same meaning as above (see after (\ref{alpha})), the function $%
A\left( x\right) $ is defined by (\ref{A}), and of course the function $%
F_{n}\left( z_{n};z_{s},s\neq n;a,b,c,d\right) $ is defined by this equation
since the expression in the curly bracket is obviously an \textit{even}
function of the $x_{s}$, hence a function of the $z_{s}=x_{s}^{2}$, for all
values of $s=1,2,...,N$.

Next, let us look at what happens \textit{infinitesimally close} to the
equilibrium configuration. To do so we set 
\begin{subequations}
\begin{equation}
x_{n}\left( t\right) =\bar{x}_{n}+\varepsilon ~\xi _{n}\left( t\right) ~,
\label{eksepsiksin}
\end{equation}%
with $\varepsilon $ infinitesimal. This of course implies the relation%
\begin{equation}
\dot{x}_{n}\left( t\right) =\varepsilon ~\dot{\xi}_{n}\left( t\right) ~,~~~%
\dot{z}_{n}\left( t\right) =2~x_{n}\left( t\right) ~\dot{x}_{n}\left(
t\right) =2~\varepsilon ~x_{n}\left( t\right) ~\dot{\xi}_{n}\left( t\right)
~.
\end{equation}

Next, let us insert the \textit{ansatz} (\ref{eksepsiksin}) in the
right-hand side of (\ref{EvEqxn}) and expand the resulting expression in
powers of $\varepsilon $. The term of zeroth order in $\varepsilon $ of
course vanishes, see (\ref{Eqxnbar}). The term of first order in $%
\varepsilon $ yields the (linearized) $t$-evolution equation 
\end{subequations}
\begin{equation}
\dot{\xi}_{n}=\mathbf{i~}\sum_{m=1}^{N}M_{nm}~\xi _{m}~,  \label{ksidot}
\end{equation}%
with the $N\times N$ matrix $\underline{M}$ defined componentwise by~(\ref%
{Mnm}) in terms of the $N$ zeros $\bar{z}_{n}=\bar{x}_{n}^{2}$ of the Wilson
polynomial $W_{N}\left( z;a,b,c,d\right) \equiv W_{N}\left(
x^{2};a,b,c,d\right) \equiv w_{2N}\left( x;a,b,c,d\right) $ (of degree $N$
in $z=x^{2}$) and of the $4$ parameters $a,b,c,d$ (or, equivalently, $\alpha
_{1},\alpha _{2},\alpha _{3},\alpha _{4}$, see (\ref{alpha})).

The terms of higher order in $\varepsilon $ yield systems of algebraic
equations satisfied by the zeros $\bar{x}_{n}$ of the Wilson polynomial $%
W_{N}\left( z\right) =w_{2N}\left( x\right) $.

The proof of \textit{Proposition~2.1} is now a consequence of the fact that
the solution of the system of linear ODEs (\ref{ksidot}) is clearly a linear
superposition (with $t$-independent coefficients) of exponentials, $\exp
\left( \mathbf{i}\mu _{m}t\right) $, where the quantities $\mu _{m}$ (with $%
m=1,2,...,N$) are the $N$ eigenvalues of the $N\times N$ matrix $\underline{M%
}$; but, due to the simultaneous validity of the relations (\ref{psicm}), (%
\ref{EvEqcm}) and (\ref{psizeros})---and to (\ref{eksepsiksin})---this
solution must also be a linear superposition (with $t$-independent
coefficients) of the $t$-dependent quantities $c_{m}\left( t\right) .$ Hence
(\ref{cmt}) implies $\mu _{m}=N\left( 2N-m+\alpha _{1}-1\right) $. \textit{%
Proposition~2.1} is thereby proven. %\end{proof}

Note that in our treatment we are implicitly assuming that the zeros $%
x_{n}\left( t\right) $ are---for all values of $t$--- all different among
themselves. This is indeed the \textit{generic} situation. Any nongeneric
event---like the \ \textquotedblleft collision" of two different zeros at
some special value of the parameter $t$---can be dealt with by appropriate
limits and in any case such possibilities---should they occur---do not
invalidate the proof of \textit{Proposition~2.1}, as reported above.

The proof of \textit{Proposition~2.2} is analogous to that of \textit{%
Proposition~2.1}; it is sketched below.

%\begin{proof}
Let $\tilde{q}_{2N}(y)=Q_{N}(y^{2}-\theta ^{2})$ be the polynomial of degree 
$2N$ in $y$ defined by~\eqref{eq:nnR8}, where $Q_{N}(z)$ is the monic
version of the Racah polynomial $R_{N}(z;\alpha ,\beta ,\gamma ,\delta )$,
see~\eqref{eq:R4}.

Suppose that the polynomial (with time-dependent coefficients) 
\begin{equation}
\tilde{\psi}_{2N}(y,t)=\tilde{q}_{2N}(y)+\sum_{m=1}^{N}c_{m}(t)\tilde{q}%
_{(2N-2m)}(y)  \label{eq:oR9}
\end{equation}%
solves the following DDE 
\begin{subequations}
\label{eq:oR8Both}
\begin{equation}
\frac{\partial }{\partial t}f(y,t)=\mathbf{i}\left[ \tilde{D}(y)(\tilde{%
\delta}_{+}-1)+\tilde{D}(-y)(\tilde{\delta}_{-}-1)-N(N+\alpha +\beta +1)%
\right] f(y,t)~,  \label{eq:oR8}
\end{equation}%
where 
\begin{equation}
\tilde{\delta}_{\pm }f(y,t)=f(y\pm 1,t)
\end{equation}%
\end{subequations}
and $\tilde{D}$ is defined by~\eqref{eq:R10}. Because $\tilde{q}_{2N}(y)$
satisfies difference equation \eqref{eq:R9}, $\tilde{q}_{2N}(y)$ is an
equilibrium solution of \eqref{eq:oR8Both}.

Using~\eqref{eq:nR8} and~\eqref{eq:nnR8}, we rewrite polynomial~%
\eqref{eq:oR9} as 

\begin{equation}
\tilde{\psi}_{2N}(y,t)=Q_{N}(\lambda
(x))+\sum_{m=1}^{N}c_{m}(t)Q_{N-m}(\lambda (x))\equiv \Psi (\lambda (x),t)~.
\label{eq:ooR9}
\end{equation}%
The last equality defines the polynomial $\Psi _{N}(z,t)$ of degree $N$
(with time-dependent coefficients). Substitution of \textit{ansatz}~%
\eqref{eq:oR9} into DDE~\eqref{eq:oR8Both} yields the system of ODEs for the
coefficients $c_{m}(t)$ 
\begin{subequations}
\begin{equation}
\dot{c}_{m}(t)=\mathbf{i}m\left( m-2N-\alpha -\beta -1\right) c_{m}(t)~,
\label{ODEc}
\end{equation}%
whose solution is 
\begin{equation}
c_{m}(t)=c_{m}(0)\exp \left[ \mathbf{i}m\left( m-2N-\alpha -\beta -1\right) t%
\right] ~.  \label{eq:oR11}
\end{equation}%
Therefore, DDE~\eqref{eq:oR8Both} indeed has a polynomial solution, with
time-dependent coefficients, given by~\eqref{eq:oR9}. Also, system of ODEs~%
\eqref{ODEc} has the unique equilibrium $c_{m}(t)\equiv 0$, which agrees
with the fact that $\tilde{q}_{2N}(y)$ is the equilibrium of DDE %
\eqref{eq:oR8Both}.

In this proof, let $z_{1}(t),z_{2}(t),\ldots ,z_{N}(t)$ denote the zeros of
the polynomial $\Psi _{N}(z,t)$ defined by~\eqref{eq:ooR9}. Then $\Psi
_{N}(z,t)=\prod_{n=1}^{N}(z-z_{n}(t))$ and 
\end{subequations}
\begin{equation}
\tilde{\psi}_{2N}(y,t)=\Psi _{N}(y^{2}-\theta ^{2})=\prod_{m=1}^{N}\left[
y^{2}-y_{m}^{2}(t)\right] ~,  \label{eq:psizeros}
\end{equation}%
where $y_{m}^{2}(t)=z_{m}(t)+\theta ^{2}$ for all $m=1,2,\ldots ,N$.

To recast DDE~\eqref{eq:oR8Both} in terms of the zeros $y_{m}(t)$ of $\tilde{\psi%
}_{2N}(y,t)$, where $m=1,2,\ldots ,N$, we compute, similarly to~\eqref{eq10}
and~\eqref{eq11}, 
\begin{equation}
\frac{\partial }{\partial t}\tilde{\psi}_{2N}(y,t)|_{y=y_{n}(t)}=-2\dot{y}%
_{n}(t)y_{n}(t)\prod_{m=1,m\neq n}^{N}\left[ y_{n}^{2}(t)-y_{m}^{2}(t)\right]
\label{eq:R15}
\end{equation}%
and 
\begin{eqnarray}
&&\left[ \tilde{D}(y)(\tilde{\delta}_{+}-1)+\tilde{D}(-y)(\tilde{\delta}%
_{-}-1)-N(N+\alpha +\beta +1)\right] \tilde{\psi}_{2N}(y,t)\Big|_{y=y_{n}(t)}
\notag \\
&=&\tilde{D}(y_{n})\tilde{\psi}_{2N}(y_{n}+1,t)+\tilde{D}(-y_{n})\tilde{\psi}%
_{2N}(y_{n}-1,t)  \notag \\
&=&\tilde{D}(y_{n})(2y_{n}+1)\prod_{m=1,m\neq n}^{N}((y_{n}+1)^{2}-y_{m}^{2})
\notag \\
&&+\tilde{D}(-y_{n})(-2y_{n}+1)\prod_{m=1,m\neq
n}^{N}((y_{n}-1)^{2}-y_{m}^{2})~.  \label{eq:R16}
\end{eqnarray}%
Taking into account relations~\eqref{eq:R15} and~\eqref{eq:R16}, we conclude
that the zeros $y_{n}(t)$ satisfy the system of ODEs 
\begin{subequations}
\label{SystODEs}
\begin{equation}
\dot{y}_{n}=\mathbf{i~}F_{n}(\underline{y})~,  \label{ODEsYn}
\end{equation}%
where $\underline{y}=(y_{1},\ldots ,y_{n})$ and 
\begin{eqnarray}
&&F_{n}(\underline{y})=-\frac{1}{2y_{n}}\Bigg\{\tilde{D}(y_{n})(2y_{n}+1)%
\prod_{\ell =1,\ell \neq n}^{N}\left( 1+\frac{1+2y_{n}}{y_{n}^{2}-y_{\ell
}^{2}}\right)   \notag \\
&&+\tilde{D}(-y_{n})(-2y_{n}+1)\prod_{\ell =1,\ell \neq n}^{N}\left( 1+\frac{%
1-2y_{n}}{y_{n}^{2}-y_{\ell }^{2}}\right) \Bigg\},~~~n=1,2,\ldots ,N~.
\label{eq:R18}
\end{eqnarray}%
This system of ODEs, (\ref{SystODEs}), can be reformulated in terms of $%
z_{n}(t)=y_{n}(t)^{2}-\theta ^{2}$ as follows: 
\end{subequations}
\begin{eqnarray}
\dot{z}_{n} &=&2y_{n}\dot{y}_{n}=-\mathbf{i}\Bigg\{\tilde{D}%
(y_{n})(2y_{n}+1)\prod_{\ell =1,\ell \neq n}^{N}\left( 1+\frac{1+2y_{n}}{%
y_{n}^{2}-y_{\ell }^{2}}\right) +[(y_{s}\rightarrow (-y_{s}))]\Bigg\}, 
\notag \\
~~~n &=&1,...,N~,  \label{Eqzndot}
\end{eqnarray}%
where the symbol $+[(y_{s}\rightarrow (-y_{s}))]$ denotes again the addition
of the expression preceding it within the curly brackets, with $y_{s}$
replaced by $(-y_{s})$ for all $s=1,2,\ldots ,N$. It is therefore obvious
that the right-hand sides of these ODEs\textbf{\ }are \textit{even}
functions of $y_{s}$, hence functions of $z_{s}=y_{s}^{2}-\theta ^{2}$, $%
s=1,2,\ldots ,N$.

Now let $\bar{z}_{1},\ldots \bar{z}_{N}$ be the $N$ zeros of the Racah
polynomial $R_{N}(z;\alpha ,\beta ,\gamma ,\delta )$, hence also the zeros
of $Q_{N}(z)$, see~\eqref{eq:R4}. Let $\pm \bar{y}_{1},\pm \bar{y}%
_{2},\ldots ,\pm \bar{y}_{N}$ be the zeros of $\tilde{q}_{2N}(y)$; these
numbers satisfy the relations $\bar{y}_{m}^{2}=\bar{z}_{m}+\theta ^{2}$. To
linearize system~\eqref{ODEsYn} about $\bar{y}_{1},\ldots ,\bar{y}_{N}$, we
set 
\begin{equation}
y_{m}(t)=\bar{y}_{m}+\varepsilon ~\eta _{m}(t)  \label{eq:ymlinearize}
\end{equation}%
with $\varepsilon $ infinitesimal. This leads to the system of linear
equations 
\begin{equation}
\dot{\underline{\eta }}=\mathbf{i~}\underline{\tilde{M}}~\underline{\eta }~,
\label{eq:R22}
\end{equation}%
where $\underline{\eta }=(\eta _{1},\ldots ,\eta _{N})$ and the matrix $%
\underline{\tilde{M}}$ is defined by 
\begin{equation}
\tilde{M}_{nm}=\frac{\partial F_{n}}{\partial y_{m}}(\underline{y})\Big|%
_{y_{k}=\bar{y}_{k},~k=1,\ldots ,N}~.  \label{eq:R23}
\end{equation}%
And it is easily seen that the explicit expressions for the components of
the matrix $\underline{\tilde{M}}$ are given by~\eqref{tildeM}.

The proof of \textit{Proposition~2.2} is now a consequence of the fact that
the solution of the system of linear ODEs \eqref{eq:R22} is clearly a linear
superposition (with $t$-independent coefficients) of exponentials, $\exp
\left( \mathbf{i}\tilde{\mu}_{m}t\right) $, where the quantities $\tilde{\mu}%
_{m}$ (with $m=1,2,...,N$) are the $N$ eigenvalues of the $N\times N$ matrix 
$\underline{\tilde{M}}$; but, due to the simultaneous validity of the
relations (\ref{eq:oR9}), (\ref{ODEc}) and (\ref{eq:psizeros})---and to (\ref%
{eq:ymlinearize})---this solution must also be a linear superposition (with $%
t$-independent coefficients) of the $t$-dependent quantities $c_{m}\left(
t\right) $, see (\ref{eq:oR9}) and (\ref{eq:oR11}). Hence (\ref{eq:oR11})
implies $\tilde{\mu}_{m}=N\left( 2N-m+\alpha _{1}-1\right) $. \textit{%
Proposition~2.1} is thereby proven. %\end{proof}

\section{Outlook}

A follow-up to the findings reported in this paper are analogous results for
the zeros of the Askey-Wilson and $q$-Racah polynomials, i. e. the "highest"
polynomials belonging to the $q$-Askey scheme \cite{KS}. Likewise, a
follow-up to the findings reported in \cite{BC2014} are analogous results
for generalized basic hypergeometric polynomials. We hope to be able to
obtain and report these results soon.

\section{Acknowledgements}

One of us (OB) would like to acknowledge with thanks the hospitality of the
Physics Department of the University of Rome ``La Sapienza" on the occasion
of three two-week visits there in June 2012, May 2013 and June-July 2014;
the results reported in this paper were obtained during the last of these
visits. The other one (FC) would like to acknowledge with thanks the
hospitality of Concordia College for a one-week visit there in November 2013.

\section{Appendix: Standard definitions and properties of the Wilson and
Racah polynomials}

In this Appendix we report for the convenience of the reader, and also to
identify the notation used throughout this paper, a number of standard
formulas associated with the Wilson and Racah polynomials. We generally
report these formulas from the standard compilations \cite{HTF1}
respectively \cite{KS} (the formulas of which are identified by the
notations (HTF1-X) respectively (KS-X), where X stands here for the notation
appropriate to identify equations in each of these compilations); in some
cases we add certain immediate consequences of these formulas which are not
explicitly displayed in these compilations nor (to the best of our
knowledge) elsewhere.

Let us recall that hereafter $\mathbf{i}$ denotes the imaginary unit, $%
\mathbf{i}^{2}=-1$.

The Pochhammer symbol is defined as follows (see (HTF1-5.1(3)): 
\begin{equation}
\left( a\right) _{0}=1~,~~~\left( a\right) _{n}=a\left( a+1\right) \cdot
\cdot \cdot \left( a+n-1\right) =\frac{\Gamma \left( a+n\right) }{\Gamma
\left( a\right) }~~~\text{for~~~}n=1,2,3,...~,  \label{Poch}
\end{equation}%
where of course $\Gamma \left( z\right) $ is the standard Gamma function.

The generalized hypergeometric function is defined as follows (see
(HTF1-5.1(2)), (KS-0.4.1)):%
\begin{equation}
_{p}F_{q}\left( \left. 
\begin{array}{c}
a_{0},a_{1},...,a_{p} \\ 
b_{1},b_{2},...,b_{q}%
\end{array}%
\right\vert z\right) =\sum_{k=0}^{\infty }\left[ \frac{\left( a_{0}\right)
_{k}\left( a_{1}\right) _{k}\cdot \cdot \cdot \left( a_{p}\right) _{k}}{%
\left( b_{1}\right) _{k}\left( b_{2}\right) _{k}\cdot \cdot \cdot \left(
b_{q}\right) _{k}}~\frac{z^{k}}{k!}\right] ~,  \label{HypFunct}
\end{equation}%
with $p$ and $q$ two positive integers (note that it is generally understood
that no one of the $q$ parameters $b_{n}$ coincide with any one of the $p+1$
parameters $a_{n}$ since in such a case these two parameters would cancel
out in the definition (\ref{HypFunct}), which would then merely reduce to
the definition of the hypergeometric function $_{p-1}F_{q-1}$). The
hypergeometric function (\ref{HypFunct}) becomes of course a polynomial in $%
z $ of degree $N$ if one of the parameters $a_{n}$ has the value $-N$ with $N
$ a positive integer---say, $a_{0}=-N$ (without loss of generality, since
the definition (\ref{HypFunct}) implies that the generalized hypergeometric
function is invariant under permutations of the $p+1$ parameters $a_{n}$ as
well as under permutations of the $q$ parameters $b_{n}$)---provided no one
of the other parameter $a_{n}$ is a negative integer smaller in modulus than 
$N$ and no one of the parameters $b_{n}$ is a negative integer (as hereafter
assumed). This is of course a simple consequence of the fact that, if $N$ is
a positive integer, $\left( -N\right) _{n}$ vanishes for $n>N,$ see (\ref%
{Poch}). Note that in this case the hypergeometric function is also a
polynomial in each of the $p$ parameters $a_{n}$ with $n=1,2,...,p.$

\subsection{Formulas for the Wilson polynomials}

The Wilson polynomial $W_{N}(z;a,b,c,d)$ with $z=x^{2}$---and its rescaled
\textquotedblleft tilded" version $\tilde{W}_{N}(z;a,b,c,d)$---are defined
as follows (see (KS-1.1.1) and (KS-1.1.4)): 
\begin{subequations}
\begin{equation}
W_{N}(z;a,b,c,d)=\left( a+b\right) _{N}~\left( a+c\right) _{N}~\left(
a+d\right) _{N}~\tilde{W}_{N}(z;a,b,c,d)~,  \label{Wtilde}
\end{equation}%
\begin{equation}
\tilde{W}_{N}(z;a,b,c,d)=_{4}F_{3}\left( \left. 
\begin{array}{c}
-N,~N+a+b+c+d-1,~a+\mathbf{i}x,~a-\mathbf{i}x \\ 
a+b,~a+c,~a+d%
\end{array}%
\right\vert 1\right) ~,  \label{Wilsona}
\end{equation}%
or equivalently (see (\ref{Poch}) and (\ref{HypFunct})),%
\begin{equation}
\tilde{W}_{N}(z;a,b,c,d)=\sum_{k=0}^{N}\left[ \frac{\left( -N\right)
_{k}~\left( N+a+b+c+d-1\right) _{k}~\left[ a;z\right] _{k}}{k!~\left(
a+b\right) _{k}\left( a+c\right) _{k}\left( a+d\right) _{k}}\right] ,
\label{Wilsonb}
\end{equation}%
where we introduced the new (modified Pochhammer) symbol 
\end{subequations}
\begin{eqnarray}
&&\left[ a;x^{2}\right] _{k}=\left( a+\mathbf{i}x\right) _{k}~\left( a-%
\mathbf{i}x\right) _{k}~,  \notag \\
&&\left[ a;z\right] _{0}=1~,  \notag \\
&&\left[ a;z\right] _{k}=\left( a^{2}+z\right) ~\left[ \left( a+1\right)
^{2}+z\right] \cdot \cdot \cdot \left[ \left( a+k-1\right) ^{2}+z\right]
,\;k=1,2,3,...~.  \notag \\
&&
\end{eqnarray}%
It is plain from this formula that $\left[ a;z\right] _{k}$ is a polynomial
of degree $k$ in $z$ (and also of degree $2k$ in $a$), hence that the Wilson
polynomial $W_{N}(z;a,b,c,d)$\textit{\ }are indeed polynomials of degree $N$
in $z$ (see (\ref{Wilsonb})).

\textit{Notational remark}. For notational simplicity we often omit to
indicate explicitly the dependence on the $4$ parameters $a$, $b$, $c$, $d$%
---provided this entails no ambiguity. $\square $

Another interesting \textit{avatar} of the Wilson
polynomials---corresponding to an alternative rescaling---reads as follows: 
\begin{subequations}
\begin{eqnarray}
&&p_{N}\left( z\right) \equiv p_{N}\left( z;a,b,c,d\right) =\frac{\left(
-1\right) ^{N}~W_{N}(z;a,b,c,d)}{\left( N+a+b+c+d-1\right) _{N}}  \notag \\
&=&\frac{\left( -1\right) ^{N}\left( a+b\right) _{N}~\left( a+c\right)
_{N}~\left( a+d\right) _{N}}{\left( N+a+b+c+d-1\right) _{N}}\cdot  \notag \\
&&\cdot \sum_{k=0}^{N}\left\{ \frac{\left( -N\right) _{k}~\left(
N+a+b+c+d-1\right) _{k}~\left[ a;z\right] _{k}}{k!~\left( a+b\right)
_{k}\left( a+c\right) _{k}\left( a+d\right) _{k}}\right\} ~.  \label{MonicW}
\end{eqnarray}%
Note that, as implied by (\ref{Wilsonb}) and (\ref{ModPoch}) (together with
the identity $\left( -N\right) _{N}=\left( -1\right) ^{N}~N!,$ see (\ref%
{Poch})), the polynomial $p_{N}\left( z\right) $ is \textit{monic}, i. e.
its highest power $z^{N}$ has \textit{unit} coefficient, implying%
\begin{equation}
\underset{z\rightarrow \infty }{\lim }\left[ \frac{p_{N}\left( z\right) }{%
z^{N}}\right] =1~.  \label{pNmonic}
\end{equation}

It is finally convenient to introduce the notation 
\end{subequations}
\begin{equation}
w_{2N}\left( x\right) \equiv w_{2N}\left( x;a,b,c,d\right)
=W_{N}(x^{2};a,b,c,d)  \label{wW}
\end{equation}%
to report the \textit{difference equation} satisfied by the Wilson
polynomials (see (KS-1.1.6), and note the minor notational changes we
introduced): 
\begin{subequations}
\label{DiffEq}
\begin{eqnarray}
&&B\left( -x\right) ~w_{2N}\left( x+\mathbf{i}\right) +B\left( x\right)
~w_{2N}\left( x-\mathbf{i}\right)  \notag \\
&=&\left[ N~\left( N+a+b+c+d-1\right) +B\left(- x\right) +B\left( x\right) %
\right] ~w_{2N}\left( x\right) ~,  \label{DiffEqa}
\end{eqnarray}%
where $B\left( x\right) \equiv B\left( x;a,b,c,d\right)$ is defined by 
\begin{equation}
B\left( x\right) \equiv B\left( x;a,b,c,d\right) =\frac{\left( a+ \mathbf{i}%
x\right) \left( b+ \mathbf{i}x\right) \left( c+ \mathbf{i}x\right) \left( d+ 
\mathbf{i}x\right) }{2\mathbf{i}x\left( 2\mathbf{i}x+ 1\right) }~.
\label{B+-}
\end{equation}%
Note that $w_{2N}\left( x\right) $ is a polynomial of degree $N$ in $x^{2},$
hence---as suggested by its notation---an \textit{even} polynomial of degree 
$2N$ in $x$.

Let us now introduce the $N$ zeros $\bar{z}_{n}\equiv \bar{z}_{n}\left(
N;a,b,c,d\right) $ of the Wilson polynomial $W_{N}\left( z\right) \equiv
W_{N}\left( z;a,b,c,d\right) $ (or $\tilde{W}_{N}\left( z\right) \equiv 
\tilde{W}_{N}\left( z;a,b,c,d\right) $ or $p_{N}\left( z\right) \equiv
p_{N}\left( z;a,b,c,d\right) $), clearly such---recalling that $p_{N}\left(
z\right) $ is monic, see (\ref{pNmonic})---that there hold the formula 
\end{subequations}
\begin{subequations}
\begin{equation}
p_{N}\left( z\right) =\prod\limits_{m=1}^{N}\left( z-\bar{z}_{m}\right) ~;
\label{DefZerosa}
\end{equation}%
as well as the $2N$ zeros, $\pm \bar{x}_{n}\equiv \pm \bar{x}_{n}\left(
N;a,b,c,d\right) =\pm \sqrt{\bar{z}_{n}},$ of the even polynomial $%
w_{2N}\left( x\right) ,$ see (\ref{wW}), of degree $2N,$ for which there
holds the analogous formula%
\begin{equation}
w_{2N}\left( x\right) =\left( -1\right) ^{N}~\left( N+a+b+c+d-1\right)
_{N}~\prod\limits_{m=1}^{N}\left( x^{2}-\bar{x}_{m}^{2}\right) ~.
\label{DefZerosb}
\end{equation}

Let us note that, for $x=\bar{x}_{n}\equiv \bar{x}_{n}\left(
N;a,b,c,d\right) $, formula (\ref{DiffEqa}) implies 
\end{subequations}
\begin{subequations}
\begin{equation}
B\left( -\bar{x}_{n}\right) ~w_{2N}\left( \bar{x}_{n}+\mathbf{i}\right)
+B\left( \bar{x}_{n}\right) ~w_{2N}\left( \bar{x}_{n}-\mathbf{i}\right) =0~,
\end{equation}%
hence, via (\ref{DefZerosb}) and (\ref{B+-}),%
\begin{equation}
\prod\limits_{m=1}^{N}\left[ \frac{\bar{x}_{n}^{2}+2\mathbf{i}\bar{x}_{n}-1-%
\bar{x}_{m}^{2}}{\bar{x}_{n}^{2}-2\mathbf{i}\bar{x}_{n}-1-\bar{x}_{m}^{2}}%
\right] =\frac{\left( 1-2\mathbf{i}\bar{x}_{n}\right) }{\left( 1+2\mathbf{i}%
\bar{x}_{n}\right) }\frac{\left( a+\mathbf{i}\bar{x}_{n}\right) \left( b+%
\mathbf{i}\bar{x}_{n}\right) \left( c+\mathbf{i}\bar{x}_{n}\right) \left( d+%
\mathbf{i}\bar{x}_{n}\right) }{\left( a-\mathbf{i}\bar{x}_{n}\right) \left(
b-\mathbf{i}\bar{x}_{n}\right) \left( c-\mathbf{i}\bar{x}_{n}\right) \left(
d-\mathbf{i}\bar{x}_{n}\right) }~,
\end{equation}%
or, equivalently,%
\begin{eqnarray}
&&\frac{\left( a+\mathbf{i}\bar{x}_{n}\right) \left( b+\mathbf{i}\bar{x}%
_{n}\right) \left( c+\mathbf{i}\bar{x}_{n}\right) \left( d+\mathbf{i}\bar{x}%
_{n}\right) }{\mathbf{i}x_{n}}\prod\limits_{\ell =1,~\ell \neq n}^{N}\left( 
\bar{x}_{n}^{2}-2\mathbf{i}\bar{x}_{n}-1-\bar{x}_{\ell }^{2}\right)  \notag
\\
&=&\frac{\left( a-\mathbf{i}\bar{x}_{n}\right) \left( b-\mathbf{i}\bar{x}%
_{n}\right) \left( c-\mathbf{i}\bar{x}_{n}\right) \left( d-\mathbf{i}\bar{x}%
_{n}\right) }{\mathbf{i}x_{n}}\prod\limits_{\ell =1,~\ell \neq n}^{N}\left( 
\bar{x}_{n}^{2}+2\mathbf{i}\bar{x}_{n}-1-\bar{x}_{\ell }^{2}\right) ~, 
\notag \\
&&
\end{eqnarray}%
or, equivalently,%
\begin{eqnarray}
\frac{\left( a-\mathbf{i}\bar{x}_{n}\right) \left( b-\mathbf{i}\bar{x}%
_{n}\right) \left( c-\mathbf{i}\bar{x}_{n}\right) \left( d-\mathbf{i}\bar{x}%
_{n}\right) }{\mathbf{i}x_{n}}\prod\limits_{\ell =1,~\ell \neq n}^{N}\left( 
\bar{x}_{n}^{2}+2\mathbf{i}\bar{x}_{n}-1-\bar{x}_{\ell }^{2}\right) && 
\notag \\
+\left[ \left( \bar{x}_{s}\rightarrow (-\bar{x}_{s})\right) \right] &=&0~, 
\notag \\
&&  \label{Eqxnbar}
\end{eqnarray}%
where (here and throughout) the symbol $+\left[ \left( \bar{x}%
_{s}\rightarrow (-\bar{x}_{s})\right) \right] $ denotes the addition of
everything that comes before it, with the replacement of $\bar{x}_{s}$ with $%
(-\bar{x}_{s})$ for all $s=1,2,\ldots ,N$.

It is plain that the left-hand sides of the $N$ algebraic equations (\ref%
{Eqxnbar}) are $N$ polynomials of degree $N+1$ in the $N$ variables $\bar{x}%
_{m}^{2}=\bar{z}_{m}$ because all terms containing an \textit{odd} power of $%
\bar{x}_{n}$ cancel out due to the addition caused by the symbol $+\left[
\left( \bar{x}_{s}\rightarrow (-\bar{x}_{s})\right) \right] $ (note
incidentally that the zeros $\bar{x}_{\ell }$ with $\ell \neq n$ only enter
as squares to begin with, see (\ref{Eqxnbar})). Hence these are $N$
polynomial equations of degree $N+1$ satisfied by the $N$ zeros $\bar{z}_{m}$
of the Wilson polynomial $W_{N}\left( z\right) $. This property of the $N$
zeros of the Wilson polynomial of degree $N$ is needed to prove our main
result (see Section 3). It obtains so easily from standard formulas for the
Wilson polynomials that it could hardly be presented as a \textit{new}
finding, although we have not been able to find it in print. The diligent
reader might wish to verify numerically that this result is correct---but
note that such a check can hardly be performed without computer assistance
(say, by Maple or Mathematica), even for quite small values of $N$.

\subsection{Formulas for the Racah polynomials}

The Racah polynomial of degree $N$ is defined by (see (KS-1.2.1) of \cite{KS}%
) 
\end{subequations}
\begin{equation}
R_{N}(\lambda (x);\alpha ,\beta ,\gamma ,\delta )=_{4}F_{3}\left( \left. 
\begin{array}{c}
-N,~N+\alpha +\beta +1,~-x,~x+\gamma +\delta +1 \\ 
\alpha +1,~\beta +\delta +1,~\gamma +1%
\end{array}%
\right\vert 1\right) ~,  \label{eq:R1}
\end{equation}%
where $N=0,1,2,\ldots ,\nu $, $\lambda (x)=x(x+\gamma +\delta +1)$ and the
parameters $\alpha ,\beta ,\gamma ,\delta $ satisfy (only) one of the
following conditions: 
\begin{equation}
\alpha +1=-\nu \mbox{ or }\beta +\delta +1=-\nu \mbox{ or }\gamma +1=-\nu ~,
\label{condRacah}
\end{equation}%
with $\nu $ a nonnegative integer. We note that $R_{N}$ is a polynomial of
degree $N$ in $\lambda (x)$: 
\begin{equation}
R_{N}(\lambda (x);\alpha ,\beta ,\gamma ,\delta )=\sum_{n=0}^{N}\frac{%
(-N)_{n}(N+\alpha +\beta +1)_{n}[\lambda (x)]_{n}}{n!(\alpha +1)_{n}(\beta
+\delta +1)_{n}(\gamma +1)_{n}}~,  \label{eq:R2}
\end{equation}%
where $(a)_{n}$ denotes the Pochhammer symbol, see \eqref{Poch}, and the
symbol $[\lambda (x)]_{n}$ is here defined as follows: 
\begin{eqnarray}
\lbrack \lambda (x)]_{n} &=&(-x)_{n}(x+\gamma +\delta +1)_{n}  \notag \\
&=&-\lambda (x)\left[ -\lambda (x)+(\gamma +\delta +1)+1\right] \left[
-\lambda (x)+2(\gamma +\delta +1)+2^{2}\right]   \notag \\
&&\cdots \left[ -\lambda (x)+(n-1)(\gamma +\delta +1)+(n-1)^{2}\right] ~.
\label{eq:R3}
\end{eqnarray}%
We define the \textit{monic} polynomial of degree $N$ in $\lambda (x)$ 
\begin{equation}
Q_{N}(\lambda (x))=\frac{(\alpha +1)_{N}(\beta +\delta +1)_{N}(\gamma +1)_{N}%
}{(N+\alpha +\beta +1)_{N}}R_{N}(\lambda (x);\alpha ,\beta ,\gamma ,\delta
)~,  \label{eq:R4}
\end{equation}%
and also introduce 
\begin{equation}
q_{2N}(x)=Q_{N}(\lambda (x))~,  \label{eq:R5}
\end{equation}%
a polynomial of degree $2N$ in $x$. By formula (KS-1.2.5), $q_{2N}(x)$
satisfies the difference equation 
\begin{subequations}
\label{eq:R8}
\begin{eqnarray}
N(N+\alpha +\beta +1)~q_{2N}(x) &=&C(x)~q_{2N}(x+1)+D(x)~q_{2N}(x-1)  \notag
\\
&&-[C(x)+D(x)]~q_{2N}(x)~,  \label{eq:R6}
\end{eqnarray}%
where 
\begin{eqnarray}
C(x) &=&\frac{(x+\alpha +1)(x+\beta +\delta +1)(x+\gamma +1)(x+\gamma
+\delta +1)}{(2x+\gamma +\delta +1)(2x+\gamma +\delta +2)}~, \\
D(x) &=&\frac{x(x-\alpha +\gamma +\delta )(x-\beta +\gamma )(x+\delta )}{%
(2x+\gamma +\delta )(2x+\gamma +\delta +1)}.  \label{eq:R7}
\end{eqnarray}

Let us now introduce the convenient parameter $\theta $ and the new variable 
$y$ as follows: 
\end{subequations}
\begin{subequations}
\begin{equation}
y=x+\theta ~,~~~\theta =\frac{\gamma +\delta +1}{2}~,  \label{ytheta}
\end{equation}
so that 
\begin{equation}
\lambda (x)=x(x+2\theta )=y^{2}-\theta ^{2}~.  \label{eq:nR8}
\end{equation}%
Because $q_{2N}(x)$ satisfies~\eqref{eq:R8}, the related polynomial 
\end{subequations}
\begin{equation}
\tilde{q}_{2N}(y)=q_{2N}(y-\theta )=Q_{N}(y^{2}-\theta ^{2})=Q_{N}(\lambda
(x))  \label{eq:nnR8}
\end{equation}%
satisfies the difference equation 
\begin{eqnarray}
N(N+\alpha +\beta +1)~\tilde{q}_{2N}(y) &=&\tilde{D}(y)~\tilde{q}_{2N}(y+1)+%
\tilde{D}(-y)~\tilde{q}_{2N}(y-1)  \notag \\
&&-\left[ \tilde{D}(y)+\tilde{D}(-y)\right] ~\tilde{q}_{2N}(y)~,
\label{eq:R9}
\end{eqnarray}%
where $\tilde{D}(y)$ is defined by~\eqref{eq:R10}.

\end{document}